\documentclass[12pt,oneside]{amsart} 
\usepackage{amssymb,amscd}
\usepackage{euscript}

      \theoremstyle{plain}
      \newtheorem{theorem}{Theorem}[section]
      \newtheorem{lemma}[theorem]{Lemma}
      
      \newtheorem{proposition}[theorem]{Proposition}

\numberwithin{equation}{section}

      \makeatletter
      \def\@setcopyright{}
      \def\serieslogo@{}
      \makeatother

\def\R{\mathbb R}
\def\Rm{\mathbb R ^m}

\def\Z{\mathbb Z}
\def\N{\mathbb N}

\def\o{\mathcal O}

\def\rm{\mathcal R^\mu}
\def\relm{\mathcal R_{\varepsilon,l}^\mu}

\newcommand{\A}{\mathcal A}
\newcommand{\Gm}{GL(m,\mathbb R)}
\newcommand{\Rmx}{\mathbb R^m_x}
\newcommand{\la}{\lambda}

\def\dist{\text{dist}}

\def\Id{\text{Id}}
\def\e{\varepsilon}
\def\a{\alpha}

\def\QED{\hfill\hfill{\square}}

\begin{document}

%\date{\today}
\author{Boris Kalinin$^\ast$}

\address{Department of Mathematics $\&$ Statistics, 
 University of South  Alabama, Mobile, AL 36688, USA}
\email{kalinin@jaguar1.usouthal.edu}

\title [Liv\v{s}ic Theorem for matrix cocycles]
{Liv\v{s}ic Theorem for matrix cocycles}

\thanks{$^{\ast}$  Supported in part by NSF grant DMS-0701292}

%%%%%%%%%%%%%%%%%%%%%%%%%%%%%%%

\begin{abstract}
We  prove the Liv\v{s}ic Theorem for arbitrary $GL(m,\R)$ cocycles.
We consider a hyperbolic dynamical system $f : X \to X$ and a H\"older 
continuous function $A: X \to GL(m,\R)$. We show that if $A$ has trivial
periodic data, i.e. $A(f^{n-1} p) ... A(fp) A(p)$ $= \Id$ for each periodic 
point $p=f^n p$, then there exists a H\"older continuous function $C: X \to GL(m,\R)$ 
satisfying $A (x)  = C(f x) C(x) ^{-1}$ for all $x \in X$. The main new ingredients
in the proof are results of independent interest on relations between the 
periodic data, Lyapunov exponents, and uniform estimates on growth of  
products along orbits  for an arbitrary H\"older function $A$.
\end{abstract}

\maketitle 

%%%%%%%%%%%%%%%%%%%%%%%%%%%%
%%%%%%%%%% Introduction %%%%%%%%%%%
%%%%%%%%%%%%%%%%%%%%%%%%%%%%

\section{Introduction}
For a hyperbolic dynamical system $f : X \to X$ and a group $G$ we consider 
the question of when a H\"older continuous function $A: X \to G$ is a {\em coboundary},
i.e. there exists a (continuous or H\"older continuous) function $C: X \to G$
satisfying 
$$
A (x)  = C(f x) C(x) ^{-1} \qquad  \text{ for all } \; x \in X.
$$
This is equivalent to the fact that the $G$-valued cocycle $\A$ generated by $A$ 
(see \eqref{gen+} and \eqref{gen-}) over the $\Z$ action generated by $f$ is 
cohomologus to the identity cocycle.
Since any coboundary $A$ must have trivial periodic data, i.e
\begin{equation}  \label{N}
\A (p,n) \overset{\text{def}}= A(f^{n-1} p) \, \cdots \, A(fp) \, A(p) = \Id \qquad 
\forall \; p \in X, \, n \in \N \; \; \text{with} \; f^n p=p , 
\end{equation} 
the question is whether this necessary condition is also sufficient.
Cocycles appear naturally in many important problems in dynamics.
A. Liv\v{s}ic was first to study cohomology of dynamical systems in 
his seminal papers \cite{Liv1,Liv2}.  In the case of Abelian $G$ 
he obtained positive answers for this and related questions.  
Similar questions for non-Abelian groups are substantially more difficult 
and, despite some progress, were not successfully resolved. 
Non-Abelian cohomology of hyperbolic systems has since been 
extensively studied, some of the highlights are 
\cite{GS,L,LMM,NP,NT95,NT98,P,PP,PW,Sch}. 
We refer the reader to \cite{LW} and to the upcoming book \cite{KN} for some 
of the most recent results and overview  of historical development in this area. 
The natural difficulty in non-Abelian Liv\v{s}ic-type 
arguments is related to the growth of the cocycle along orbits. 
In particular, the sufficiency of condition \eqref{N} was established
when $G$ is compact or when $A$ is either sufficiently close to identity 
or satisfies some growth assumptions. For example, specific 
{\em localization assumptions} are given in \cite{LW} for various cases 
of groups and metrics on them. 
\vskip.2cm

In this paper we prove the sufficiency of \eqref{N} for an arbitrary $GL(m,\R)$
cocycle, which has been a long standing open problem. We also obtain an
important result for cocycles with uniformly bounded periodic data. Our 
theorems cover most classes of groups with interesting applications, except 
for groups of diffeomorphisms. To prove these theorems we establish new 
relations between the periodic data, Lyapunov exponents, and uniform 
estimates of the growth for an arbitrary H\"older cocycle. These results 
are of independent interest and have wide applicability.

To include various classes of hyperbolic systems  $f : X \to X$ and streamline 
the notations we formulate explicitly the property that we will use.

\vskip.2cm

\noindent{\bf Definition.} We call orbit segments $x, fx, ... , f^n x$ and $p, fp, ... , f^n p$
{\em exponentially $\delta$ close with exponent $\la >0$} if for every $i=0, ... , n$ we have
\begin{equation}  \label{d-close-traj}
\dist (f^i x, f^i p) \le \delta \cdot \exp \left( -\la\, \min\{i,n-i \} \right) .
\end{equation}

\vskip.2cm

\noindent{\bf Definition.} We say that a homeomorphism $f$ of a 
metric space $X$ satisfies {\em closing property} if there exist 
$c \, , \la, \delta_0 >0$ such that for any $x \in X$ and $n>0$ with 
$\dist (x, f^n x) < \delta_0$ there exists a point $p \in X$ with 
$f^n p =p$ such that the orbit segments $x, fx, ... , f^n x$ and 
$p, fp, ... , f^n p$ are exponentially $\delta = c \,  \dist (x, f^n x)$ close 
with exponent $\la$ and there exists a point $y \in X$ 
such that for every $i = 0, ..., n$
\begin{equation}  \label{mp}
\dist (f^i p, f^i y)  \le  \delta \, e^{-\la i} \quad \text{and} \quad
\dist (f^i y, f^i x)  \le \delta \, e^{-\la (n-i)}.
\end{equation}

\vskip.2cm

Anosov Closing Lemma and the local product structure yield the closing 
property for smooth hyperbolic systems such as hyperbolic automorphisms 
of tori and nilmanifolds, Anosov diffeomorphisms, and locally maximal 
hyperbolic sets (basic sets of axiom-A systems) \cite{KH}. Another class 
satisfying the closing property includes symbolic dynamical systems such 
as subshifts of finite type. 
\vskip.2cm

We now state our main result, the Liv\v{s}ic Theorem for matrix cocycles.
Recall that a homeomorphism is called {\em topologically transitive} 
if it has a dense orbit.

\begin{theorem} \label{livsic} Let $f$ be a topologically transitive 
homeomorphism
of a compact metric space $X$ satisfying the closing property. Let
$A: X \to GL(m,\R)$ be an $\a$-H\"older function such that 
$$
A(f^{n-1} p) \, \dots \, A(fp) \, A(p) = \Id \qquad 
\forall \; p \in X, \, n \in \N \; \; \text{with} \; f^n p=p.
$$
Then there exists an $\a$-H\"older function $C : X \to GL(m,\R)$  
such that 
\begin{equation}  \label{C}
A (x)  = C(f x) C(x) ^{-1} \qquad  \text{ for all } \; x \in X.
\end{equation}

\end{theorem}

\vskip.2cm

\noindent{\bf Remark.} 
Note that a value of $C$ at a point $x$ uniquely determines by \eqref{C}
the values of $C$ on the orbit of $x$. 
Hence, by the topological transitivity of $f$, $C$ is unique up to a translation, 
i.e. any other $C'$ satisfying \eqref{C} is of the form $C'(x)=C(x)B$ for 
some $B \in GL(m,\R)$. Also,  \cite[Theorem 2.4]{NT98} 
implies  that such $C$ is smooth if so are $A$ and $(X,f)$. 

\vskip.2cm

\noindent{\bf Remark.} As we note in the end of the proof, if $A $ takes values 
in a closed subgroup $G$ of $GL(m,\R)$ then $C$ can be naturally chosen to take
values in $G$. Thus Theorem \ref{livsic} holds if $GL(m,\R)$ is replaced by such a
group $G$. In fact, the theorem holds for any connected Lie group $G$ as follows
 from the remark after the next theorem.

\vskip.2cm

Next we consider a more general case when the periodic data is not 
trivial but is uniformly bounded, for example is contained in a compact 
subgroup. In this case we prove that the cocycle itself is also bounded.

\begin{theorem} \label{bound} Let $f$ be a transitive homeomorphism
of a compact metric space $X$ satisfying the closing property and let 
$A: X \to GL(m,\R)$ be an $\a$-H\"older function. Suppose that there exists
a compact set $K \subset GL(m,\R)$ such that $\A(p,n) \in K$
for all $p \in X$ and $n \in \N$ with $f^n p=p$. Then there exists a compact set $K'$ such that $\A(x,n) \in K'$ 
for all $x \in X$ and $n \in \Z$.
\end{theorem}

In particular, this theorem allows one to obtain further cohomology information 
for $GL(m,\R)$ cocycles with uniformly bounded periodic data by using 
results obtained  in \cite{Sch} for cocycles that distort a distance on the 
group in a bounded fashion. 
For further results on cocycles with bounded or conformal periodic data see subsequent paper \cite{KaS}.
\vskip.2cm

\noindent{\bf Remark.} 
For a cocycle with values in a connected Lie group $G$ Theorem \ref{bound} 
can be applied to the adjoint representation. For example, if the periodic data 
is trivial \eqref{N} then the theorem  implies that  all $Ad \, (\A(x,n))$ are 
uniformly bounded and hence 
the cocycle distorts a right invariant metric on $G$ in a bounded fashion. 
It follows from  \cite{Sch} or classical arguments \cite{Liv2},
 \cite[Theorem 5.3.1]{KN} that Theorem \ref{livsic} holds for such $G$.

\vskip.2cm

To prove Theorems \ref{livsic} and \ref{bound} we first establish the following 
growth estimates for a cocycle in terms of its periodic data. This result gives 
new tools for further study of cohomology for 
non-Abelian cocycles, in particular for the case when the periodic data 
has exponents close to zero. We think that Theorem \ref{slow}  will also 
be useful for various problems in smooth dynamics of hyperbolic systems 
and actions, such as existence of invariant geometric structures and rigidity.

\begin{theorem} \label{slow} Let $f$ be a homeomorphism of a compact metric space 
$X$ satisfying the closing property and let $\A$ be a H\"older $GL(m,\R)$ cocycle over $f$.  
Let $\, \chi_{min}$ and $\chi_{max}$ be real numbers such that for every periodic point $p$
every eigenvalue $\rho$ of $\A (p,n)$ satisfies 
$ \chi_{min} \le \frac 1n \log |\rho | \le  \chi_{max}$, 
where $n$ is the period of $p$. Then for any $\e > 0$ there exists a constant $c_\e$ 
such that for all $x \in X$ and $n \in \N$
\begin{equation}  \label{estA}
 \| \mathcal A(x,n) \| \le c_{\e} \exp(n \chi_{max} +\e n) \quad \text{and} \quad
 \| \mathcal A(x,n)^{-1} \| \le c_{\e} \exp(-n \chi_{min} +\e n) .
\end{equation}
\end{theorem}

\vskip.2cm

The proof of this theorem relies on our next result which resembles 
Theorem 3.1 in \cite{SW} on approximation of Lyapunov exponents of a 
hyperbolic invariant measure for a diffeomorphism that follows earlier 
results in \cite{K}.  Note that in our case 
there is no assumption on hyperbolicity of the cocycle and, in fact, 
our main application is to cocycles with all Lyapunov exponents equal to zero.

\begin{theorem} \label{periodic} Let $f$ be a homeomorphism
of a compact metric space $X$ satisfying the closing property, 
let $\A$ be a H\"older $GL(m,\R)$ cocycle over $f$, and let $\mu$ be  
an ergodic invariant measure for $f$. Then the Lyapunov exponents
 $\chi_1 \le ... \le \chi_m$ (listed with multiplicities) of $\A$ with respect 
 to $\mu$ can be approximated 
by the Lyapunov exponents of $\A$ at periodic points. More precisely,
for any $\e >0$ there exists a periodic point $p \in X$ for which the Lyapunov 
exponents $\chi_1^{(p)} \le ... \le \chi_m^{(p)}$ of $\A$ satisfy
$|\chi_i-\chi_i^{(p)}|<\e$ for $i=1, \dots , m$.
\end{theorem}
 \vskip.2cm

\noindent{\bf Remark.} Theorems \ref{slow} and \ref{periodic} use 
only a weaker version of the closing property without the existence of a point $y$.
Also, $\delta = c \, \dist (x, f^n x)$ in the closing property could be replaced
by $\delta = c \, \dist (x, f^n x)^\beta$ with $\beta >0$. The proofs of Theorems \ref{bound}, \ref{slow}, and 
\ref{periodic} work in the same way with proper modifications of  exponents. 
Similarly, Theorem \ref{livsic} holds in this case with $C$ being $(\a \beta)$-H\"older.
\vskip.2cm

\noindent{\bf Remark.} More generally, Theorems \ref{livsic}, \ref{bound}, \ref{slow}, 
and \ref{periodic} hold for an extension $\A$ of $f$ by linear transformations
of a vector bundle $\mathcal B$ over $X$. The arguments are essentially
identical since we compare the values of $A$ and related structures only 
at nearby points. This can be done if one can identify fibers at nearby points
H\"older continuously via local trivialization or connection. In particular, the 
theorems apply to the derivative cocycle of a smooth hyperbolic system, as
well as to its restriction to a H\"older continuous invariant distribution,
without any global trivialization assumptions.

\vskip.2cm

We would like to thank Victoria Sadovskaya, Ralf Spatzier, and  Anatole Katok 
for helpful comments and suggestions.

%%%%%%%%%%%%%%%%%%%%%%%%%%%%%
%%%%%%%%%%      Preliminaries   %%%%%%%%%%
%%%%%%%%%%%%%%%%%%%%%%%%%%%%%

\section{Cocycles over $\Z$ actions}\label{preliminaries}

In this section we review some basic definitions and facts of the 
Oseledec theory of cocycles over $\Z$ actions. 
We use \cite{BP} as a general refernce.

\subsection{Cocycles}

Let $f$ be an invertible transformation of a  space $X$. 
A function $\mathcal A: X \times \Z \to\Gm$ is called a {\em linear cocycle} or 
a {\em matrix--valued cocycle} over $f$ if for all $x \in X$ and $n,k \in \Z$  
we have $\mathcal A(x,0)=\Id$ and
\begin{equation}  \label{cocycle}
\mathcal A(x,n+k)=\mathcal A(f^k x, n ) \cdot \mathcal A(x,k).
\end{equation}
We consider only matrix-valued cocycles and simply call them {\em cocycles}.
Any cocycle $\mathcal A(x,n)$ is uniquely determined by its generator $A : X \to \Gm$, 
which we sometimes also call cocycle. The generator is defined by $A(x) = \mathcal A(x,1)$, 
and the cocycle can be reconstructed from its generator as follows, for any $n>0$
\begin{equation}  \label{gen+}
\A(x,n)= A(f^{n-1}x) \cdots  A(fx) \cdot A(x),
\end{equation}
\begin{equation}  \label{gen-}
\A(x,-n)= A(f^{-n}x)^{-1} \cdots  A(f^{-2}x)^{-1} \cdot A(f^{-1}x)^{-1}
= \A(f^{-n} x,n)^{-1}.
\end{equation}

\vskip.2cm

A cocycle $\A$ over a homeomorphism $f$ of a metric space $X$ is called 
$\a$-H\"older if its generator $A : X \to GL(m,\R)$ is H\"older continuous
with exponent $\a$. To consider this notion we need to introduce a metric 
on $GL(m,\R)$, for
example as follows
\begin{equation}  \label{distG} 
\dist_{GL(m,\R)} (A, B) = \| A  - B \|  + \| A^{-1}  - B^{-1} \| , \quad \text {where}
\end{equation}
$$
\| A \| = \sup \{ \| Au \| \cdot \| u\| ^{-1} : \; 0 \not= u \in \Rm \} \, .
$$
We note that on any compact set in $GL(m,\R)$ the norms  $\| A^{-1}  \|$ and $\| B^{-1} \|$ 
are uniformly bounded and hence this distance is Lipschitz 
equivalent to $\| A  - B \|$. 
Therefore, for a compact $X$, a cocycle $\A$ is $\a$-H\"older if and only 
if $\| A (x)  - A(y) \| \le c \, \dist (x,y)^\a$ for all $x,y \in X$. For a non-compact
$X$ certain caution is needed as in the proof of Theorems \ref{livsic} and \ref{bound}.

\vskip.2cm

\subsection{Lyapunov exponents and Lyapunov metric}\label{SLyapunovmetric} 

Cocycles can be considered in various categories. Even though in this paper we 
mostly study H\"older cocycles, a general theory is developed for measurable 
cocycles over measure preserving transformations. 

\begin{theorem} [Oseledec Multiplicative Ergodic Theorem, see \cite{BP} Theorem 3.4.3] 
\label{MET}  
Let $f$ be an invertible ergodic measure-preserving transformation of a Lebesgue 
probability measure space $(X,\mu)$. Let $\mathcal A$ be a measurable cocycle 
whose generator satisfies $\log\|A(x)\|\in L^1(X,\mu)$ and $\log\|A(x)^{-1}\|\in L^1(X,\mu)$. 
Then there exist numbers $\chi_1 < \dots < \chi_l$, an $f$-invariant set $\rm$ with 
$\mu (\rm)=1$, and an $A$-invariant Lyapunov decomposition of $\Rm$ for $x\in \rm$ 
$$
\Rmx= E_{\chi_1}(x)\oplus\dots\oplus E_{\chi_l}(x)
$$ 
with $\dim E_{\chi_i}(x)=m_i$, such that for any $i=1,...,l$ and any $0 \not= v\in E_{\chi_i}(x)$ one has 
$$ \lim_{n\to{\pm \infty}} n^{-1} \log\|\mathcal A(x,n) v \|=  \chi_i   \quad \text{and} \quad
 \lim_{n\to{\pm \infty}} n^{-1}  \log\det\mathcal A(x,n)= \sum_{i=1}^lm_i  \chi_i. 
$$
\end{theorem}

\noindent{\bf Definitions.}  The numbers $\chi_1,\dots,\chi_l$
are called the {\it Lyapunov exponents} of $A$ and the
dimension $m_i$ of the space $E_{\chi_i}(x)$ is called the {\it
multiplicity} of the exponent $\chi_i$.  The points of the set $\rm$
are called {\it regular}. 

\vskip.2cm

We denote the standard scalar product in $\Rm$ by $<\cdot,\cdot>$. For
a fixed $\e >0$ and a regular point $x$ we introduce {\it the $\e$-Lyapunov 
scalar product (or metric)}  $<\cdot,\cdot>_{x,\e}$ in $\Rm$ as follows.
For $u\in E_{\chi_i}(x),\,\,v\in E_{\chi_j}(x),\,\,i\neq j$  we set  $<u,v>_{x,\e}\, =0$.
For $i=1,\dots,l$ and $u,v\in E_{\chi_i}(x)$  we define 
$$
<u,v>_{x,\e} \, =m \sum_{n\in\Z}<\mathcal A(x,n)u,\mathcal
A(x,n)v>\exp(-2\chi_i n -\e |n|). 
$$ 
Note that the series converges exponentially for any regular $x$. The constant $m$ in front of
the conventional formula is introduced for more convenient comparison with 
the standard scalar product. 
Usually,  $\e$ will be fixed and we will denote $<\cdot,\cdot>_{x,\e}$
simply by $<\cdot,\cdot>_x$ and call it the {\it Lyapunov scalar product}. The norm 
generated by this scalar product is called the {\em Lyapunov norm}
and is denoted by $\|\cdot\|_{x,\e}$ or $\|\cdot\|_x$.

\vskip.2cm

We summarize below some important properties of the Lyapunov scalar product 
and norm, for more details see  \cite[Sections 3.5.1-3.5.3]{BP}. A direct calculation 
shows \cite[Theorem 3.5.5]{BP} that for any regular $x$ and any $u\in E_{\chi_i}(x)$ 
\begin{equation}  \label{estAEi}
\exp(n \chi_i -\e |n|) \|u\|_{x,\e} \le 
\| \mathcal A(x,n)u\|_{f^n x,\e} \le
\exp(n \chi_i+\e |n|)\|u\|_{x,\e} \quad \forall n \in \Z,
\end{equation}
\begin{equation}  \label{estAnorm}
\exp(n \chi -\e n) \le \| \mathcal A(x,n) \|_{f^n x \leftarrow x} \le \exp(n \chi+\e n) 
\quad \forall n \in \N,
\end{equation}
where $\chi = \chi _l$ is the maximal Lyapunov exponent and 
$\| \cdot \|_{f^n x \leftarrow x}$ is the operator norm with respect to the
Lyapunov norms. It is defined for any matrix $A$ and any regular points $x,y$ as follows
$$
\| A \| _{y \leftarrow x}=\sup \{ \| Au \|_{y,\e} \cdot \| u\|_{x,\e} ^{-1} : \; 0 \not= u \in \Rm \}.
$$

We emphasize that, for any given $\e>0$, Lyapunov scalar product and 
Lyapunov norm are defined only for regular points with 
respect to the given measure. They depend only measurably on the 
point even if the cocycle is H\"older. Therefore, comparison with the standard
norm becomes important. The uniform lower bound follows easily from the 
definition: $\|u\|_{x,\e}\ge \|u\|$. The upper bound is not uniform,
but it changes slowly along the regular orbits 
\cite[Proposition 3.5.8]{BP}: there exists a 
measurable function $K_\e (x)$ defined on the set of regular points $\rm$ 
such that 
\begin{equation}  \label{estLnorm}
\| u \| \le  \| u \|_{x,\e} \le K_\e(x) \|u\| \qquad \forall x \in \rm, \; \; 
\forall u \in \Rm \qquad \text{and} 
\end{equation}
\begin{equation}  \label{estK}
 K_\e(x) e^{-\e |n|}  \le K_\e(f^n x) \le  K_\e(x) e^{\e |n|}  \qquad 
  \forall x \in \rm, \; \forall  n \in \Z.
\end{equation}
These estimates are obtained in \cite{BP} using the fact that 
$\|u\|_{x,\e}$ is {\em tempered}, but they can also be checked directly using the definition of $\|u\|_{x,\e}$ on each Lyapunov space and noting that angles between the spaces change slowly. 

For any matrix $A$ and any regular points $x,y$ inequalities 
\eqref{estLnorm} and \eqref{estK} yield 
\begin{equation}  \label{estMnorm}
K_\e (x) ^{-1} \| A \| \le \| A \|_{y \leftarrow x} \le K_\e (y) \| A \| \, .
\end{equation}

When $\e$ is fixed we will usually omit it and write $K(x)=K_\e (x)$.
For any $l>1$ we also define the following sets of regular points
\begin{equation}  \label{Pset}
\relm = \{x \in \rm : \; \; K_\e(x) \le l \}.
\end{equation}
Note that $\mu (\relm) \to 1$ as $l \to \infty$.
Without loss of generality we can assume that the set $\relm$ is compact
and that Lyapunov splitting and Lyapunov scalar product are continuous on
$\relm$. Indeed, by Luzin theorem we can always find a subset of $\relm$ 
satisfying these properties with arbitrarily small loss of measure (in fact, for 
standard Pesin sets these properties are automatically satisfied).

%%%%%%%%%%%%%%%%%%%%%%%%%%%%%%
%%%%%%%%%%%%%   Proofs   %%%%%%%%%%%%

%%%%%%%%%%%%%%%%%%%%%%%%%%%%%%%%%%%%

\section{Proof of Theorem \ref{periodic}}  \label{proof-periodic}

We begin with Lemma \ref{mainest} below which gives a general estimate of 
the norm of $\A$ along {\em any} orbit segment close to a regular one.
In fact, its proof does not use the measure $\mu$ and relies only on the 
estimates for $\A$ and the Lyapunov norm along the orbit segment 
$x, fx, ... , f^n x$ that follow from the fact that $x, f^n x \in \relm$.

\begin{lemma} \label{mainest}
Let $\A$ be an $\a$-H\"older cocycle over a homeomorphism $f$ of a compact
metric space $X$ and let $\mu$ be an ergodic measure for $f$ with the largest 
Lyapunov exponent $\chi$. Then for any positive $\lambda$ and $\e$
satisfying $\la > \e /\a$ there 
exists  $c>0$ such that 
for any $n \in \N$, any regular point $x$ with both $x$ and $f^n x$ in $\relm$,
and any point $y \in X$ such that the orbit segments $x, fx, ... , f^n x$ and 
$y, fy, ... , f^n y$ are exponentially $\delta$ close with exponent $\lambda$ we have
\begin{equation}  \label{d-close-coc}
\| \A(y,n) \|_{f^n x \leftarrow x} \le e^{c \,  l \delta^\a} e^{n (\chi+\e) } \le  
e^{2n \e + c \, l \delta^\a}  \, \| \A(x,n) \|_{f^n x \leftarrow x}
\qquad  \text{and}
\end{equation}
\begin{equation}  \label{d-close-norm}
\| \A(y,n) \| \le l \, e^{c \,  l \delta^\a} e^{n (\chi+\e) }
\le l^2 e^{2n \e + c \,  l \delta^\a} \, \| \A(x,n) \| . \quad
\end{equation}
The constant $c$ depends only on the cocycle $\A$ and on the 
number $(\a \la - \e)$.
\end{lemma}

\vskip.2cm

\proof
We denote $x_i=f^i x$ and $y_i=f^i y$, $i=0, ... , n$, and
estimate the Lyapunov norm
$$
\| \A(y,n) \|_{x_n \leftarrow x_0} = \| A(y_{n-1}) \; ... \; A(y_1) \, A(y_0) \|_{x_n \leftarrow x_0}=
$$
$$
= \| A(x_{n-1}) \, [A(x_{n-1}) ^{-1} A(y_{n-1}) ] \; ... \; A(x_0)\, [A(x_0)^{-1} A(y_0) ] \|_{x_n \leftarrow x_0} \le
$$
$$
\| A(x_{n-1}) \|_{x_n \leftarrow x_{n-1}}  \|A(x_{n-1}) ^{-1} A(y_{n-1}) \|_{x_{n-1} \leftarrow x_{n-1}}  ... 
\| A(x_0) \|_{x_1 \leftarrow x_0}  \| A(x_0)^{-1} A(y_0)  \|_{x_0 \leftarrow x_0} .
$$
\vskip.2cm

Since $\| A(x_i) \|_{x_{i+1} \leftarrow x_i} \le e^{\chi+\e}$ by \eqref{estAnorm}, 
where $\chi$ is the maximal exponent of $\A$ at $x$, we conclude that
\begin{equation} \label{mainest1}
\| \A(y,n) \|_{x_n \leftarrow x_0} \le  e^{n (\chi + \e)} \,
\prod_{i=0}^{n-1}  \| A(x_i)^{-1} A(y_i)  \|_{x_i \leftarrow x_i} 
\end{equation}   
To estimate the product term we consider $D_i =  A(x_i)^{-1} A(y_i) - \Id$.
Since $A(x)$ is $\a$-H\"older on the compact space $X$, and hence $\| A(x)^{-1} \|$ 
is uniformly bounded, we obtain using the closeness of the orbit segments that
\begin{equation} \label{mainestDi}
\| D_i \| \le \| A(x_i)^{-1} \| \cdot  \| A(y_i) - A(x_i) \| 
\le c' \dist (x_i, y_i) ^ \a \le c' \left( \delta e^{-\la \min \{i, n-i\} } \right)^ \a ,
\end{equation}
where the constant $c'$ depends only on the cocycle $\A$.
Since both $x$ and $f^n x$ are in $\relm$ we have $K(x_i) \le l e^{\e \min\{i,n-i \}}$
by \eqref{estK} and \eqref{Pset}. Hence for the Lyapunov norms we can conclude that
\begin{equation} \label{DestL}
\| D_i \|_{x_i \leftarrow x_i} \le  K(x_i) \|  D_i \| \le l e^{\e \min\{i,n-i \}} \,   \| D_i \|
\le l e^{\e \min\{i,n-i \}} \, c' \delta^ \a e^{-\la \a \min \{i, n-i\} } 
\end{equation}   
\begin{equation} \label{Rest}
\text{and} \quad \| A(x_i)^{-1} A(y_i) \|_{x_i \leftarrow x_i} \le 1 + \| D_i \|_{x_i \leftarrow x_i} \le
1 +  c' l  \,  \delta^\a \, e^{(\e -\a \la) \, \min\{i,n-i \} } . 
\end{equation}   
Now using \eqref{mainest1} and \eqref{Rest} we obtain
$$
\log (\| \A(y,n) \|_{x_n \leftarrow x_0}) - n (\chi + \e) \le 
\sum _{i=0}^{n-1} \log \| A(x_i)^{-1} A(y_i) ] \|_{x_i \leftarrow x_i} \le
$$
$$
\le c' l \delta^\a \sum _{i=0}^{n-1}  \exp \left[(\e -\a \la) \, \min\{i,n-i \} \right]
\le c \,  l \delta^\a
$$
since the sum is uniformly bounded due to the assumption $\e < \a \la$. 
The constant $c$ depends only on the cocycle $\A$ and on $(\a \la - \e)$.
We conclude using \eqref{estAnorm} that 
\begin{equation} \label{mainestL}
\| \A(y,n) \|_{x_n \leftarrow x_0} \le e^{c \,  l \delta^\a} e^{n (\chi + \e)} \le  
e^{2n \e + c \,  l \delta^\a}  \, \| \A(x,n) \|_{x_n \leftarrow x_0} .
\end{equation} 
Since $K(x_0) \le l$ and $K(x_n) \le l$ we can also estimate the standard norm 
$$
 \| \A(y,n) \| \le K(x_0) \| \A(y,n) \|_{x_n \leftarrow x_0} \le 
l e^{c \,  l \delta^\a} e^{n (\chi + \e)} \le 
l e^{2n \e + c \,  l \delta^\a} \, \| \A(x,n) \|_{x_n \leftarrow x_0} \le
$$
\begin{equation} \label{mainestE}
\le l e^{2n \e + c \,  l \delta^\a} \, K(x_n) \, \| \A(x,n) \|
\le l^2 e^{2n \e + c \,  l \delta^\a} \, \| \A(x,n) \| .
\end{equation}   
Estimates \eqref{mainestL} and \eqref{mainestE} complete 
the proof of Lemma \ref{mainest}. 

$\QED$
\vskip.2cm

The main part of the proof of Theorem \ref{periodic} is the following proposition
which gives approximation for the largest Lyapunov exponent of $\A$. 
We use it to complete the proof of Theorem \ref{periodic} at the end of 
Section \ref{proof-periodic}.

Let $f$ be a homeomorphism of a compact metric space $X$ satisfying 
the closing property with exponent $\la$, let $\A$ be an $\a$-H\"older $GL(m,\R)$ 
cocycle over $f$, and let $\mu$ be an ergodic invariant measure for $f$. 
We denote by $\chi$ the largest Lyapunov exponent of $\A$ with respect 
to $\mu$. Similarly, for any periodic point $p$ we denote by $\chi^{(p)}$ the 
largest Lyapunov exponent of $\A$ at $p$. We set
$\e _0 = \min \{ \la \a , (\chi - \nu)/2) \}$, where $\nu < \chi$ is the second 
largest Lyapunov exponent with respect to $\mu$.
In the case when $\chi$ is the only Lyapunov exponent of $\A$ with respect 
to $\mu$, we take $\e _0 =  \la \a$.

\begin{proposition} \label{largestexp}
Let $f$, $\A$, $\mu$, and $\e_0$ be as above. 
Then for any positive $l$ and $\e < \e _0$ there exist $N,\delta >0$ such that if a
periodic orbit $p, fp, ... , f^n p=p$ is exponentially $\delta$ close to an
orbit segment $x, fx, ... , f^n x$, with $x, f^n x$ in $\relm$ and $n>N$, 
then $|\chi-\chi^{(p)}| \le 3\e$. 
\end{proposition}
\vskip.2cm

\proof 
To estimate $\chi^{(p)}$ from above we apply Lemma \ref{mainest} with $p=y$. 
Note that the largest exponent at $p$ satisfies
$$
\chi^{(p)} \le n^{-1} \log \| \A(p,n) \| .$$
From the first inequality in \eqref{d-close-norm} we obtain that 
$$
n^{-1} \log \| \A(p,n) \|  \le \chi + \e + n^{-1} \log ( l \, e^{c \,  l \delta^\a})
.$$ 
We conclude that $\chi^{(p)} \le \chi + 2\e $ provided that $\delta$ is small
enough and  $n$ is large enough compared to $l$.

To estimate $\chi^{(p)}$ from below we will estimate the growth of 
vectors in a certain cone $K \subset \Rm$ invariant under $\A(p,n)$.
As in Lemma \ref{mainest} we first consider an arbitrary orbit segment
close to a regular one. Let $x$ be a point in  $\relm$ and $y \in X$ 
be a point such that the orbit segments $x, fx, ... , f^n x$ and $y, fy, ... , f^n y$ 
are exponentially $\delta$ close with exponent $\la$. We denote $x_i=f^i x$ and 
$y_i=f^i y$, $i=0, ... , n$. For each $i$ we have orthogonal splitting $\Rm = E_i \oplus F_i$, 
where $E_i$ is the Lyapunov space at $x_i$ corresponding to the largest Lyapunov 
exponent $\chi$ and $F_i$ is the direct sum of all other Lyapunov spaces at $x_i$
corresponding to the Lyapunov exponents less than $\chi$. 
 For any vector $u \in \Rm$ we denote by 
$u = u' + u^\perp$ the corresponding splitting with $u' \in E_i$ and $u^\perp \in F_i$, 
the choice of $i$ will be clear from the context. To simplify notations, we
write $\| . \| _i$ for the Lyapunov norm at $x_i$.
For each $i=0, ... , n \,$ we consider cones
$$
K_i = \{ u \in \Rm : \| u^\perp \| _i \le  \| u' \| _i\} \quad \text{and} \quad
K_i^\eta = \{ u \in \Rm : \| u^\perp \| _i \le (1-\eta)  \| u' \| _i\}
$$ 
with $\eta >0$.
We will consider the case when $\chi$ is {\em not} the only Lyapunov exponent 
of $\A$ with respect to $\mu$. Otherwise $F_i=\{ 0 \}$, 
$K_i^\eta = K_i =\Rm$, and the argument becomes simpler.
 Recall that $\e < \e _0 = \min \{ \la \a , (\chi - \nu)/2) \}$, where $\nu < \chi$ 
is the second largest Lyapunov exponent of $\A$ with respect to $\mu$.

\begin{lemma} \label{cone}
In the notations above, for any regular set $\relm$ there exist $\eta, \delta > 0$ 
such that if $x, f^n x \in \relm$ and the orbit segments 
$x, fx, ... , f^n x$ and $y, fy, ... , f^n y$ are exponentially $\delta$ close with exponent 
$\la$ then for every  $i=0, ... , n-1$ we have $A(y_i) (K_i ) \subset K_{i+1}^\eta$ 
and  $\| \left( A(y_i) u \right) ' \|_{i+1} \ge e^{\chi -2\e} \| u' \| _i \,$ for any $u \in K_i$.
\end{lemma}

\proof 
We fix $0 \le i < n$ and write
$$
A(y_i) = A(y_i) A(x_i)^{-1} A(x_i)  = (\Id + D_i ) \, A(x_i),
$$
where similarly to \eqref{mainestDi} we have
\begin{equation} \label{DE}
\| D_i \| = \| A(y_i) A(x_i)^{-1} - \Id \|  \le  \| A(y_i) - A(x_i) \| \, \| A(x_i)^{-1} \|
\le c_1 \dist (x_i, y_i) ^ \a .
\end{equation}   
For any  $u = u' + u^\perp \in K_i$ we consider $v= A(x_i) u \,$ and its splitting
$v = v' + v^\perp$ with $v' \in E_{i+1}$ and $v^\perp \in F_{i+1}$. Then by \eqref{estAEi}
we have $\| v \| _{i+1} \le e^{\chi + \e} \| u \| _i$ as well as
$$
 \| v' \| _{i+1} = \| A(x_i) u' \| _{i+1} \ge e^{\chi - \e} \| u' \| _i \quad \text{and} \quad
  \| v^\perp \| _{i+1} = \| A(x_i) u^\perp \| _{i+1} \le e^{\nu + \e} \| u^\perp \| _i \, .
$$
Now we consider $w = A(y_i) u = (\Id + D_i ) v = v + D_i v$ and 
its splitting $w = w' + w^\perp$ with $w' \in E_{i+1}$ and $w^\perp \in F_{i+1}$.
Then we have  
\begin{equation} \label{wv}
w' = v' + (D_i v)' \quad \text{and} \quad
 w^\perp = v^\perp + (D_i v)^\perp.
\end{equation}   
 Now using \eqref{DE}  we obtain 
$$
 \| D_i v \| _{i+1} \le \| D_i \|_{x_{i+1} \leftarrow x_{i+1}} \| v \| _{i+1} \le 
 K (x_{i+1}) \| D_i \| \,  e^{\chi + \e} \| u \| _i \le  
$$
$$
 l e^{\e \min \{i+1, n-i-1\} } \, c_1 \dist (x_i, y_i) ^ \a  \;  e^{\chi + \e} \, \sqrt{2} \, \| u' \| _i  \; ,
$$
as both $x_0$ and $x_n$ are in $\relm$.  Since 
$\dist (x_i, y_i) \le \delta e^{-\la \min \{i, n-i\} }$ we conclude that
\begin{equation} \label{Dv}
 \| D_i v \| _{i+1} \le 
 \sqrt{2} \, l c_1 e^\e \delta^\a e^{(-\la \a + \e) \min \{i, n-i\} }   \| u' \| _i 
 \le c_2 l \, \delta^\a   \| u' \| _i  \; ,
\end{equation} 
since $-\la \a + \e <0$.  Now using \eqref{wv} and \eqref{Dv} we obtain that
for small enough $\delta$
$$
 \| w' \| _{i+1} \ge 
e^{\chi - \e} \| u' \| _i - c_2 l \, \delta^\a   \| u' \| _i 
\ge e^{\chi - 2\e}   \| u' \| _i \, ,
$$
which gives the inequality in the lemma. Similarly we obtain an upper estimate 
\begin{equation} \label{w'}
 \| w' \| _{i+1} \le 
e^{\chi + \e} \| u' \| _i + c_2 l \, \delta^\a   \| u' \| _i  \le c_3  \| u' \| _i \, .
\end{equation}   
Finally, from \eqref{wv} we have
$$
\| w' \| _{i+1} \ge \| v' \| _{i+1} -  \| D_i v \| _{i+1} \quad \text{and} \quad
  \| w^\perp \| _{i+1} \le \| v^\perp \| _{i+1} +  \| D_i v \| _{i+1} \; ,
$$
so that using \eqref{Dv} again we can estimate 
$$
 \| w' \| _{i+1} -  \| w^\perp \| _{i+1} \ge 
 \| v' \| _{i+1} - \| v^\perp \| _{i+1} - 2  \| D_i v \| _{i+1} \ge
$$
$$
\ge e^{\chi - \e} \| u' \| _i - e^{\nu + \e} \| u^\perp \| _i -2 c_2 l \, \delta^\a   \| u' \| _i
\ge (e^{\chi - \e} - e^{\nu + \e} -2 c_2 l \, \delta^\a )  \| u' \| _i \ge \eta' \, \| u' \| _i
$$
for any fixed $\eta' < (e^{\chi - \e} - e^{\nu + \e})$ provided that $\delta$ is small enough. 
Now using \eqref{w'} we conclude that 
$ \| w' \| _{i+1} -  \| w^\perp \| _{i+1} \ge  \eta \, \| w' \| _{i+1}$
with $\eta =  \eta'/c_3$. This shows that $w \in K_{i+1}^\eta$ and hence 
$A(y_i) (K_i ) \subset K_{i+1}^\eta$.
This completes the proof of Lemma \ref{cone}.

$\QED$

\vskip.2cm

We now apply this lemma to the periodic orbit $p, fp, ... , f^n p=p$ and
conclude that $A(p,n) (K_0 ) \subset K_{n}^\eta$. 
Since the Lyapunov splitting and Lyapunov metric are continuous on the
compact set $\relm$, the cones $K^\eta_0$ and $K^\eta_n$ are close
if $x$ and $f^n x$ are close enough. Therefore we can ensure that
$K_{n}^\eta  \subset K_0$ if $\delta$ small enough and thus 
$A(p,n) (K) \subset K$ for $K=K_0$. Finally, using the norm estimate
in the lemma we obtain for any $u \in K$  
$$
\| A(p,n)\, u \| _{n} \ge \| (A(p,n)\, u)' \| _{n} \ge  e^{n(\chi - 2\e)} \| u' \| _0 
\ge \frac1{\sqrt{2}} e^{n(\chi - 2\e)} \| u \| _0 \ge \frac12 e^{n(\chi - 2\e)} \| u \| _n
$$
since Lyapunov norms at $x$ and $f^n x$ are close if $\delta$ is small enough.
Since $A(p,n)\, u \in K$ for any $u \in K$, we can iteratively apply $A(p,n)$ and
use the inequality above to estimate the largest Lyapunov exponent at $p$
$$
\chi^{(p)} \ge \chi (u) = \lim _{k \to \infty} \frac1{kn} \log \| \A(p,kn) u \|_{n} \ge \frac1{n} 
\lim _{k \to \infty} \frac1{k} \log  \left( \left(\frac12 e^{n(\chi - 2\e)}\right)^k \| u \| _n \right) \ge
$$
$$
\ge  \frac1{n} \left[n(\chi - 2\e) - \log 2 \right] + \frac1{n} \lim _{k \to \infty} \frac{\| u \| _n}{k}
\ge (\chi - 2\e) -  \frac{\log 2}n \ge \chi - 3\e
$$
provided that $n$ is large enough. This gives the desired lower estimate and
completes the proof of Proposition \ref{largestexp}. 

$\QED$

\vskip.1cm

We will now complete the proof of Theorem \ref{periodic}.
We apply Proposition \ref{largestexp} to cocycles $\wedge ^i \, \A$ induced by 
$\A$ on the $i$-fold exterior powers $\wedge ^i \, \Rm$, for $i= 1, ... , m$.
This trick is related to Ragunatan's proof of Multiplicative Ergodic Theorem 
\cite[Section 3.4.4]{BP} and was also used in \cite{SW}. We note that the 
largest Lyapunov exponent of $\wedge ^i \, \A$ is equal to 
$(\chi _m + ... + \chi_{m-i+1})$, where $\chi _1 \le ... \le \chi_{m}$ are the
Lyapunov exponents of $\A$ listed with multiplicities.

For any positive $\e < \e_0$ we choose $l$ so that $\mu (R) >0$, where $R$ 
is the intersection of the sets $\relm$ for all cocycles $\wedge ^i \, \A$, $i= 1, ... , m$. 
We may assume that $\mu$ is not atomic since the theorem is trivial otherwise.
We take $x\in R$ to be a non-periodic 
point with $\mu (B_r (x) \cap R) >0$ for any $r>0$, where $B_r (x)$ is the ball of
radius $r$ centered at $x$.  Then by Poincar\'{e} recurrence there exist iterates $f^n x$,
with $n$ growing to infinity, returning to $R$ arbitrarily close to $x$. Therefore, 
by the closing property, for any $\delta >0$ there exists a periodic point $p$ with 
$f^n p =p$ such that orbit segments $x, fx, ... , f^n x$ and $p, fp, ... , f^n p$ 
are exponentially $\delta$ close with exponent $\la$.  
Then Proposition \ref{largestexp} implies that for small enough  $\delta $
such a periodic point $p$ gives the approximation  
$$|(\chi _m + ... + \chi_{m-i+1})-(\chi^{(p)} _m + ... + \chi^{(p)}_{m-i+1})| \le 3\e$$
 for all $i= 1, ... , m$.  This yields the 
simultaneous approximation for all $\chi_i$, $i= 1, ... , m$, and completes
the proof of Theorem \ref{periodic}.

$\QED$

%%%%%%%%%%%%%%%%%%%%%%%%%%%%%%%%%%%%%%%%%%

\section{Proof of Theorem \ref{slow}}  \label{proof-slow}

The assumption on the eigenvalues of $\A (p,n)$ implies that all
Lyapunov exponents of $\A$ at all periodic orbits are in the interval
$[\chi_{min}, \chi_{max}]$. It follows from Theorem \ref{periodic}  that the Lyapunov exponents 
of $\A$ are in $[\chi_{min}, \chi_{max}]$ for {\em any} ergodic $f$-invariant 
measure. Such control on exponents gives the desired uniform estimates
on the growth of the norm of the cocycle. This uses a result on subadditive
sequences obtained in \cite{S}. We formulate here a weaker version
sufficient for our purposes, which appeared with a short proof in \cite{RH}.
 
\vskip.2cm

 \cite[Proposition 3.4]{RH}
{\it Let $f : X \to X$ be a continuous map of a compact metric space.
 Let $a_n : X \to \R$, $n \geq 0$, be a sequence of continuous functions 
such that 
\begin{equation} \label{3.1}
    a_{n+k} (x) \le a_n (f^k (x)) + a_k (x) 
    \;\text{ for every }x \in X,\;\; n, k \geq 0
\end{equation}    
 and such that there is a sequence of continuous functions $b_n : X \to \R$, $n \geq 0$,
satisfying 
\begin{equation} \label{3.2}
    a_n (x) \le a_n (f^k (x)) + a_k (x) + b_k (f^n (x))
     \;\text{ for every } x \in X,\;\; n, k \geq 0.  
\end{equation}   
If  $\;\inf _n \left( \frac1n \int _X a_n d \mu \right) < 0\;$ 
for every ergodic $f$-invariant measure, then there
is $N \geq 0$ such that 
$a_N (x) < 0$ for every $x \in X$.}
\vskip.2cm

We take $\e >0$ and apply this result to $\;a_n(x)=\log \| \A(x,n)\| - (\chi_{max} + \e) n$. 
It is easy to see that $a_n$ satisfy  \eqref{3.1}. Then the Subadditive Ergodic 
Theorem (or \cite[Theorem 3.5.5]{BP}, or equations \eqref{estAnorm},\eqref{estK}, and \eqref{estMnorm}) 
implies that for every $f$-invariant ergodic measure $\mu$, its maximal
exponent $\chi$, and $\mu$-a.e. $x \in X$ 
$$
\inf _n \,\frac1n \int _X a_n d \mu =\, \lim _{n\to \infty} \,\frac 1n {a_n(x)}    
= \chi - (\chi_{max} + \e) <0 \; ,
$$
and thus the assumptions on $a_n$ are satisfied. Taking into account \eqref{3.1}
we see that \eqref{3.2} holds once $a_n(x) \le a_{n+k}(x) + b_k(f^nx)$
is satisfied. This is easily verified for $b_k (x)=\log \| \A(x,k) ^{-1} \|$ since
by the cocycle identity \eqref{cocycle} we have
$$
\| \A(x,n) \| \le \| \A(f^n x,k)^{-1} \| \cdot \| \A(x,n+k) \|.
$$
We conclude from the proposition above that for any $\e >0$ there exists $N_\e$ 
such that $a_{N_\e}(x)<0$, i.e. $\| \A(x,N_\e) \| \le e^{(\chi_{max} + \e) N_\e}$ 
for all $x \in X$. Hence \eqref{estA} is satisfied for all $x$ in 
$X$ and $n$ in $\N$, where $c_\e = \max \| \A(x,k) \|$ with the maximum 
taken over all $x \in X$ and $1 \le k < N_\e$. 
The other estimate in \eqref{estA} is obtained similarly, for example by
applying the same argument to the cocycle generated by $A^{-1}$ over $f^{-1}$.
This completes the proof of Theorem \ref{slow}.
$\QED$

%%%%%%%%%%%%%%%%%%%%%%%%%%%%%%%%%

\section{Proof of Theorem \ref{livsic} and Theorem \ref{bound}}  \label{proof-livsic}

We follow the usual approach of extension along a dense orbit.  
Our proof is similar to the one in \cite{LW} with some 
modifications for the case of bounded periodic data. 
The main difference is that Theorem \ref{slow} enables us
to apply the following proposition. This allows us to complete 
the proof without extra assumptions on the cocycle $\A$. 

\begin{proposition} \label{closeA} 
Let $f$ be a homeomorphism of a compact metric space $X$ and let $A$ be 
an $\a$-H\"older $GL(m,\R)$ cocycle over $f$ such that for some $\e > 0$ and $c_\e$ 
\begin{equation}  \label{egrow}
\|  \A (x,n) \|  \le c_{\e} e^{\e n} \quad \text{ and }  \quad
 \|  \A(x,n)^{-1} \| \le c_{\e} e^{\e n} 
\qquad \forall \; x \in X \; ,  n \in \N.
\end{equation} 
Then for any $\la > 2\e /\a$ there exists a constant $c$, which depends only on $A$, $c_\e$, 
and $(\a \la - 2\e)$, such that for any $\delta$ and any orbit segments $x, fx, ... , f^n x$
and $y, fy, ... , f^n y$
\begin{equation}  \label{Axy}
\text{if} \; \; \dist (f^i x, f^i y)  \le  \delta e^{-\la i}, \; i = 0, ..., n, \; \;  
\text{ then } \; \| \A (x,n)^{-1}  \A (y,n) - \Id \,\| \leq c \, \delta^\a
\end{equation}
$$
\text{and if} \; \; \dist (f^i x, f^i y)  \le    \delta e^{-\la (n-i)}, \; i = 0, ..., n, \; \;  \text{then} \; \| \A (x,n)  \A (y,n)^{-1} - \Id \,\| \leq c \, \delta^\a.
$$
 \end{proposition}
 
\proof

We will consider the case when $\dist (f^i x, f^i y)  \le   \delta e^{-\la i}$
for $i = 0, ..., n$. The other case can be proved similarly.
Denoting $D_i = A(f^{i}x)^{-1} \, A(f^{i}y) - \Id$, $i = 0, ..., n-1$, we can write
$$
\A (x,n)^{-1}   \A (y,n) =  
\A (x,n-1)^{-1} \, A(f^{n-1}x)^{-1} \, A(f^{n-1}y) \, \A (y,n-1)  =
$$
$$
=  \A (x,n-1)^{-1}  (\Id + D_{n-1}) \, \A (y,n-1)  =
$$
$$
=  \A (x,n-1)^{-1}   \A (y,n-1) +   \A (x,n-1)^{-1}  D_{n-1} \, \A (y,n-1)   =
$$
$$
= ... =  \Id + \sum_{i=0}^{n-1} \A (x,i)^{-1}  D_i \, \A (y,i) \, .
$$
Therefore using the assumption \eqref{egrow} we obtain
$$
 \| \A (x,n)^{-1}   \A (y,n) - \Id \,\| \le 
 \sum_{i=0}^{n-1} \| \A (x,i)^{-1} \| \cdot \| D_i \| \cdot \| \A (y,i) \| \le 
 \sum_{i=0}^{n-1} (c_{\e} e^{\e i})^2  \, \| D_i \| \, .
$$ 
Similarly to \eqref{mainestDi} we can estimate
$$
\| D_i \| = \| A(f^{i}x)^{-1} \, A(f^{i}y) - \Id \, \| \le 
  c_1 \dist(f^{i}x,f^{i}y)^\a \le c_1 \delta^\a e^{-\a \la i} \, .
$$
Using the two estimates above and the assumption $\la > 2\e /\a$ we conclude that 
$$
 \| \A (x,n)^{-1}   \A (y,n) - \Id \,\| \ \le 
 \sum_{i=0}^{n-1} c_1 \, c_{\e}^2 \, \delta^\a \, e^{(2\e -\a \la) i} \le c \, \delta^\a \, ,
$$ 
where the constant $c$ depends only on $A$, $c_\e$, and $(\a \la - 2\e) >0$.

$\QED$
\vskip.2cm

We will now prove Theorems \ref{bound} and \ref{livsic}. 
Note that the condition on the periodic data of $\A$ in either theorem implies that 
the assumptions of Theorem \ref{slow} are satisfied with $\chi_{min}=\chi_{max}=0$ and 
hence \eqref{estA} gives \eqref{egrow} with any $\e >0$. Therefore, we can
take $\e < \a \la /2$, where $\la$ is the exponent in the closing property for $f$.

In the proof we will abbreviate $d_G = \dist _{GL(m,\R)}$.
Since $f$ is transitive, there exists a point $z\in X$ with dense 
orbit $\o =\{ f^k z \}_{k \in \Z}$. We will show that $d_G (\A (z,k), \Id)$
is uniformly bounded in $k \in \Z$. Since $\o$ is dense and $\A$ is
continuous this implies that $d_G (\A (x,n), \Id)$ is uniformly bounded 
in $x \in X$ and $n \in \Z$. This yields Theorem \ref{bound}.
\vskip.1cm

Consider any two points of $\o$ for which $\dist(f^{k_1}z, f^{k_2}z)< \delta_0$,
where $\delta_0$ is as in the closing property. 
Assume $k_1<k_2$ and denote $x=f^{k_1}z$ and $n=k_2-k_1$, so that 
$\delta =\dist(x, f^n x)< \delta_0$. 
By the closing property there exist points $p,y \in X$ 
with $f^n p =p$ such that for $i = 0, ..., n$
$$
\dist (f^i y, f^i p)  \le  c \, \delta \, e^{-\la i} \qquad \text {and} \qquad
\dist (f^i y, f^i x)  \le  c \, \delta \, e^{-\la (n-i)} \, .
$$
Now using Proposition \ref{closeA}  we obtain  
\begin{equation}  \label{P-Id} 
 \| \A (p,n)^{-1}   \A (y,n) - \Id \,\| \leq c_1 \delta^\a \quad \text {and} \quad
\|   \A (x,n) \A (y,n)^{-1}  - \Id \,\| \leq c_1 \delta^\a .
\end{equation}
We want to show that these inequalities imply that there exists $c_2$ such that
\begin{equation}  \label{dpyx} 
d_G ( \A (p,n), \A (y,n))\leq c_2 \delta^\a \quad \text{and} \quad d_G ( \A (y,n), \A (x,n)) \leq c_2 \delta^\a
\end{equation}
uniformly in $x,p,y,n$. We use the following simple estimate.

\begin{lemma} \label{matrix}
If $d_G(A, \Id \, ) \le M$ and either $\| A^{-1}B - \Id \, \| \le \xi $ or 
$\| AB^{-1} - \Id \, \| \le \xi $, with $\xi < 1/2$,  then $d_G(A, B) \le 3(M+1) \xi$.
\end{lemma}

\proof
We prove the first case, the second case follows similarly. 
From the assumption we have $\| A \| \le M+1$ and $\| A^{-1} \| \le M+1$. Then
$$
\| A-B \| \le \| A \| \cdot \| \Id -A^{-1} B \| \le (M+1) \xi .
$$
Denoting $Y=\Id - A^{-1}B$ we obtain
$B^{-1}A= (\Id - Y)^{-1} = \Id+ Y + Y^2 + ...\,$. Then
$$
\| B^{-1}A - \Id \, \| \le \sum_{k=1}^{\infty} \| Y ^k \| 
\le \sum_{k=1}^{\infty} \xi ^k = \frac{\xi}{1-\xi}\le 2\xi \qquad \text{and}
$$
$$
\| A^{-1}-B^{-1} \| \le \| A^{-1} \| \cdot \| \Id -B^{-1} A \| \le (M+1)2 \xi \, , 
$$
so that $d_G(A, B) =\| A-B \|+ \| A^{-1}-B^{-1} \| \le 3(M+1) \xi$. 

$\QED$
\vskip.2cm

Since  the periodic data is in a compact subset of $GL(m,\R)$ there exists $c_0$
so that
\begin{equation}  \label{PB} 
d_G ( \A (p,n), \Id) \leq c_0
\end{equation}
for all $p$ and $n$. Now Lemma \ref{matrix} and the first equation in \eqref{P-Id}
give the first equation in \eqref{dpyx} which implies, in particular, that 
$d_G ( \A (y,n), \Id)$ is also uniformly bounded. Then the lemma and the second 
equation in \eqref{P-Id} give the second equation in \eqref{dpyx}.
This establishes \eqref{dpyx}, which implies that 
\begin{equation}  \label{dpx} 
d_G ( \A (p,n), \A (x,n))\leq 2c_2 \delta^\a \qquad \text{and hence} 
\end{equation}
\begin{equation}  \label{dpx'} 
d_G (\A (x,n), \Id) \le c_0 + 2c_2 \delta^\a \le c_3 \hskip1.8cm
\end{equation}
for all $x \in \o$ and $n \in \Z$ with $\delta = \dist (x, f^n x) < \delta_0 $. 
The case of negative $n$ follows from the corresponding estimate for positive $n$.

By density of $\o$ we can take its finite piece $\o_L =\{ f^k z \}_{k \in [-L,L]}$ which 
forms a $\delta_0$ net in $X$ and choose $c_4 = \max_{k \in [-L,L]}d_G( \A(z,k) , \Id)$. 
Then for any $N \in \Z$ there exists $k \in [-L,L]$ such that $\dist(f^{k}z, f^{N}z)< \delta_0$. 
Denoting  $x=f^{k}z$ and $n=N-k$ we have $\dist(x, f^n x)< \delta_0$, so that
\eqref{dpx'} applies. The cocycle property \eqref{cocycle} gives
$$
\A(z,N) = \A (x,n) \, \A (z,k)  .
$$
Since the distance from $\Id$ to the terms on the right is bounded by $c_3$ and $c_4$ 
we conclude that $d_G (\A(z,N), \Id)$ is also uniformly bounded. This completes the proof
of Theorem \ref{bound}.
\vskip.2cm

To prove Theorem \ref{livsic} we define a function $C: \o \to GL(m,\R)$ by 
$C(f^n z) = \A(z,n)$. Note that $C$ satisfies \eqref{C} for $x \in \o$ and that 
$d_G (C, \Id)$ is uniformly bounded by the previous argument.  
It remains to show that $C$ is $\a$-H\"older on $\o$ 
with uniform constant and hence extends uniquely to an $\a$-H\"older function on $X$, 
which also satisfies \eqref{C}. Indeed, consider any $x\in \o$  and $n \in \Z$ with 
$\dist(x, f^n x)=\delta < \delta_0$. Since $\A (p,n) = \Id$ by the assumption, 
using \eqref{dpx} we obtain 
$$
\| C(f^n x) C(x) ^{-1} - \Id \, \| < d_G(C(f^n x) C(x) ^{-1}, \Id) = 
 d_G(\A (x, n), \Id) \leq 2c_2 \delta^\a  .
$$
Now, since $d_G (C, \Id)$ is uniformly bounded, Lemma \ref{matrix} gives the desired 
H\"older continuity of $C : \o \to GL(m,\R)$. This completes the proof of Theorem \ref{livsic}.

Note that if the function $A : X \to GL(m,\R)$ takes values 
in a subgroup $G \subset GL(m,\R)$
then so does the function $C$ on $\o$ and, if $G$ is closed, so does
the extension  $C : X \to GL(m,\R)$.
$\QED$

%%%%%%%%%%%%%%%%%%%%%%%%%%%%%%%%%%%%
%%%%%%%%               bibliography                %%%%%%%%%

\end{document}